\documentclass[12pt, oneside, reqno]{amsart}
\usepackage{hyperref}   	
\usepackage{geometry}
\numberwithin{equation}{section}                		
\usepackage{graphicx}				
\usepackage{amsmath,amssymb,amsthm,amscd}
\usepackage{ascmac}
\usepackage{indentfirst}
\usepackage{bm}
\usepackage{color}
\usepackage{longtable}
\usepackage{cases}
\usepackage[all]{xy}
\usepackage{url}
\usepackage[OT2,T1]{fontenc}
\usepackage{stackengine,scalerel}
\DeclareSymbolFont{cyrletters}{OT2}{wncyr}{m}{n}
\DeclareMathSymbol{\Sha}{\mathalpha}{cyrletters}{"58}
\DeclareMathOperator{\Gal}{\mathrm{Gal}}

\newtheorem{thm}{\textbf{Theorem}}[section]
\newtheorem*{thm*}{Theorem}
\newtheorem{mainthm}{Theorem}

\newtheorem{cor}[thm]{\textbf{Corollary}}
\newtheorem*{cor*}{\textbf{Corollary}}
\newtheorem{conj}[thm]{\textbf{Conjecture}}
\newtheorem*{conj*}{\textbf{Conjecture}}
\newtheorem{defn}[thm]{\textbf{Definition}}
\newtheorem{prop}[thm]{\textbf{Proposition}}
\newtheorem*{prop*}{\textbf{Proposition}}
\newtheorem{lem}[thm]{\textbf{Lemma}}

\theoremstyle{definition}
\newtheorem{example}[thm]{\textbf{Example}}
\newtheorem{rem}[thm]{Remark}

\renewcommand{\proofname}{\textit{Proof.}}

\newcommand{\al}[1][]{\alpha}
\newcommand{\be}[1][]{\beta}
\newcommand{\om}[1][]{\omega}
\newcommand{\Q}[1][]{\mathbb{Q}}
\newcommand{\Z}[1][]{\mathbb{Z}}
\newcommand{\R}[1][]{\mathbb{R}}
\newcommand{\C}[1][]{\mathbb{C}}
\newcommand{\F}[1][]{\mathbb{F}}
\newcommand{\Sub}[1][]{\mathcal{S}}
\newcommand{\hh}[1][]{\mathcal{H}}
\newcommand{\D}[1][]{\mathbf{D}_{\mathrm{crys}}}
\newcommand{\LL}[1][]{\mathcal{L}}
\newcommand{\HH}[1][]{H_{\mathrm{Iw}}}
\newcommand{\pr}[1][]{\mathrm{pr}}
\newcommand{\ra}[1][]{r_\mathrm{an}}
\newcommand{\ral}[1][]{r_\mathrm{alg}}
\newcommand{\alg}[1][]{\mathrm{alg}}
\newcommand{\Sel}[1][]{\mathrm{Sel}}
\newcommand{\g}[1][]{\gamma}
\newcommand{\Ker}[1][]{\mathrm{Ker}}
\newcommand{\Res}[1][]{\mathrm{Res}}
\newcommand{\semi}[1][]{\mathrm{ss}}
\newcommand{\Hom}[1][]{\mathrm{Hom}}
\newcommand{\loc}[1][]{\mathrm{Loc}}
\newcommand{\cyc}[1][]{\mathrm{cyc}}
\newcommand{\Reg}[1][]{\mathrm{Reg}}
\newcommand{\Fr}[1][]{\mathrm{Frob}}
\newcommand{\cl}[1][]{\mathrm{Cl}}
\newcommand{\ur}[1][]{\mathrm{ur}}
\newcommand{\A}[1][]{\mathcal{A}}
\newcommand{\z}[1][]{\mathbf{z}}
\newcommand{\zz}[1][]{\mathbf{z}_{\Q}}
\newcommand{\mcO}[1][]{\mathcal{O}}
\newcommand{\tens}[1]{%
\mathbin{\mathop{\otimes}\limits_{#1}}%
}
\title[]{On $p$-adic $L$-functions of elliptic curves and the ideal class groups of the division fields}

\author{Naoto Dainobu}
\address{Faculty Of Mathematics, Kyushu University, Motooka 744, Nishi-ku Fukuoka 819-0395, Japan}
\email{dainobu.naoto.819@m.kyushu-u.ac.jp}
\begin{document}
\maketitle
\begin{abstract}
Let $E$ be an elliptic curve defined over $\Q$ and $F$ be $\Q$ or an imaginary quadratic field with certain conditions. In this article, we study the ideal class group $\cl(F_E)$ of the $p$-division field $F_E:=F(E[p])$ of $E$ over $F$ for an odd prime number $p$. More precisely, we investigate the non-vanishing of the $E[p]$-component in the semi-simplification of $\cl(F_E)/p\cl(F_E)$ as an $\F_p[\Gal(F_E/F)]$-module when $E[p]$ is an irreducible $\Gal(F_E/F)$-module. When the analytic rank of $E$ over $F$ is $1$, we establish a new relationship between the non-vanishing of the $E[p]$-component and the $p$-divisibility of a certain $p$-adic analytic quantity associated with $E$. The quantity is defined by the leading coefficient of the cyclotomic $p$-adic $L$-function of $E$ when $F=\Q$ and by that of Bertolini--Darmon--Prasanna's anticyclotomic $p$-adic $L$-function of $E$ when $F$ is the imaginary quadratic field.
\end{abstract}

\section{Introduction and main results}\label{introresults}

\subsection{Introduction}\label{intro}
One of the central themes in number theory is the study of relationships between algebraic objects such as ideal class groups and analytic objects such as $L$-functions. The famous result of Herbrand and Ribet \cite{Her32,Rib76} reveals a certain relationship between ideal class groups of cyclotomic fields and special values of the Riemann zeta function. They considered the character decompositions of the ideal class groups as Galois modules and investigated the non-vanishing of each component in terms of the Riemann zeta values. 

In the unpublished note \cite{P17}, Prasad considered a non-abelian analogue of Herbrand and Ribet's result for division fields of elliptic curves. He tried to establish a relationship between ideal class groups of division fields and special values of the Hasse--Weil $L$-functions of elliptic curves. We can partially obtain such a relationship as an immediate consequence of his subsequent work \cite{PS21} with Shekhar under the Birch and Swinnerton-Dyer (BSD) conjecture.

In this paper, for elliptic curves of analytic rank $1$, we establish a new relationship between the ideal class groups of their division fields and certain {\it $p$-adic analytic} quantities associated with the elliptic curves. Our new relationship partially refines the previous works \cite{P17,PS21}. 

We now explain this in a little more detail. Let $E$ be an elliptic curve over $\Q$. We fix a minimal Weierstrass model of $E$ over $\Q$. Let $E[p]$ be the group of $p$-torsion points of $E$. We consider the $p$-division field $\Q_E:=\Q(E[p])$ of $E$, which is generated over $\Q$ by the coordinates of all the points in $E[p]$. We can think of $\Q_E$ as an analogue of cyclotomic fields in Herbrand and Ribet's result since cyclotomic fields are generated by torsion points of the multiplicative group scheme $\mathbb{G}_m$.

Hereafter, we assume the irreducibility of $E[p]$ as a $\Gal(\Q_E/\Q)$-module, which is almost always satisfied. Under this assumption, the extension $\Q_E/\Q$ is non-abelian. 

Let $\cl(\Q_E)$ be the ideal class group of $\Q_E$. We put $H:=\cl(\Q_E)/p\cl(\Q_E)$ and consider its semi-simplification 
\begin{align}\label{introsemi}
H^{\semi} = \bigoplus_{M\hspace{1mm} :\hspace{1mm} \text{irr}} M^{\oplus r(M)}
\end{align}
as an $\F_p[\Gal(\Q_E/\Q)]$-module. Here $M$ in the above direct product runs through all irreducible $\F_p[\Gal(\Q_E/\Q)]$-modules and the non-negative integer $r(M)$ denotes the multiplicity of the $M$-component in $H^{\semi}$. Our interest is in an investigation of the non-vanishing of the $E[p]$-component, or the non-negative integer $r(E[p])$ in (\ref{introsemi}) in terms of analytic quantities associated with $E$. 

Let $L(E, s)$ be the Hasse--Weil $L$-function associated to $E$ over $\Q$, which has an analytic continuation to the whole complex plane by the modularity theorem. We write $L^{\ast}(E, 1)$ for the leading coefficient of the Taylor expansion of $L(E, s)$ at $s=1$, and $\ra(E/\Q)$ for the vanishing order $\mathrm{ord}_{s=1}L(E, s)$. We consider the algebraic part 
\[
\displaystyle L^{\ast}(E, 1)_{\alg}:=\frac{L^{\ast}(E, 1)}{\Omega_E \mathrm{Reg}(E/\Q)}
\]
of $L(E, 1)$, where $\Omega_E$ is the N\'eron real period of $E$ and $\mathrm{Reg}(E/\Q)$ is the regulator of $E$ over $\Q$. It is conjectured that $L^{\ast}(E, 1)_{\alg} \in \Q$ and this is known to be true when $\ra(E/\Q) \leq 1$.

In \cite{P17}, Prasad expected the existence of a certain relationship between the non-vanishing $r(E[p])\neq 0$ and the $p$-divisibility of $L^{\ast}(E, 1)_{\alg}$. In \cite{PS21}, he and Shekhar studied the non-vanishing of $r(E[p])$ in terms of the $p$-Selmer group $\Sel(\Q, E[p])$ of $E$. Their result \cite[Theorem 3.1]{PS21} yields the following implication under the $p$-part of the BSD conjecture, which realizes a part of Prasad's expectation in \cite{P17} : 
\begin{align}\label{intromot}
v_p(L^{\ast}(E, 1)_{\alg})\geq 1 \hspace{2mm}\Rightarrow \hspace{2mm} r(E[p]) \neq 0.
\end{align}
Here we also assume $L^{\ast}(E, 1)_{\alg} \in \Q$ and write $v_p$ for the $p$-adic valuation which is normalized as $v_p(p)=1$.

One may expect that the converse of (\ref{intromot}) is also true, comparing with Herbrand and Ribet's result. However, we can observe many counterexamples to the converse of (\ref{intromot}) at least when $\ra(E/\Q) = 1$ as we explain below. In addition, actually, the right-hand side of (\ref{intromot}) always holds when $\mathrm{rank}_{\Z}E(\Q) \geq 2$ by \cite[Theorem 3.1]{PS21}. However, it seems unreasonable to expect that $v_p(L^{\ast}(E, 1)_{\alg}) \geq 1$ always holds in such a situation. Therefore, there seems to be a large gap between the two statements in (\ref{intromot}) in general.

\begin{example}\label{introcounterexam}
Let $p=13$ and $E$ be the elliptic curve with LMFDB \cite{lmfdb} label 43.a1, which is defined by the equation $y^2+y=x^3+x^2$ and satisfies $\ra(E/\Q) = 1$.

We can check $13 \nmid L^{\prime}(E, 1)_{\alg}$ by numerical computation. On the other hand, the author showed that $r(E[13]) \neq 0$ in \cite[Example 5.2]{Da23}. This gives a counterexample to the converse of (\ref{intromot}). As a side note, the $13$-part of the BSD conjecture for $E$ holds in this case. 
\end{example}

In the above example, we have a rational point $P \in E(\Q)$ of infinite order with certain special conditions. In \cite{Da23}, we call such a rational point an {\it everywhere unramified} rational point and show that such a point yields a non-trivial element in the $E[p]$-component of $H^{\semi}$. Therefore, the large contributions of everywhere unramified rational points of $E$ are overlooked in (\ref{intromot}). In our first main result (Theorem \ref{intromain1}), when $\ra(E/\Q)=1$, we add such missing contributions to the left-hand side of (\ref{intromot}), using certain quantities defined by the cyclotomic $p$-adic $L$-functions of $E$. 

In the second main result (Theorem \ref{intromain2}), we also treat the ideal class group of the $p$-division field $K(E[p])$ of $E$ over an imaginary quadratic field $K$. As a result, we obtain a similar result to Theorem \ref{intromain1}, using the anticyclotomic $p$-adic $L$-function of $E$ over $K$.

\subsection{Main results}
Here we state our main results. Let $E$ be an elliptic curve over $\Q$ which has good reduction at an odd prime number $p$. We fix a minimal Weierstrass model of $E/\Q$. 

Let $F$ be $\Q$ or an imaginary quadratic field. We write $F_E:=F(E[p])$ for the $p$-division field of $E$ over $F$ and $\cl(F_E)$ for the ideal class group of $F_E$. In the following, we assume that $E[p]$ is an irreducible $\F_p[\Gal(F_E/F)]$-module. We consider the semi-simplification $H_F^{\semi}$ of $H_F:=\cl(F_E)/p\cl(F_E)$ as a $\Gal(F_E/F)$-module as in (\ref{introsemi}), and write $r_F(E[p])$ for the multiplicity of the $E[p]$-component in $H_F^{\semi}$. When $F=\Q$, we write $r_{\Q}(E[p])$ as $r(E[p])$ for short.

In the following, we fix an algebraic closure $\bar{\Q}$ of $\Q$. We also fix an algebraic closure $\bar{\Q}_p$ of $\Q_p$, and fix an embedding $\bar{\Q}\hookrightarrow \bar{\Q}_p$.
\subsubsection{Result on division fields over $\Q$}
First, we consider the case where $F=\Q$. Let $\al, \be \in \bar{\Q}_p$ be the two roots of the Hecke polynomial $X^2-a_p(E)X+p=0$ with $v_p(\al) \leq v_p(\be)$. Here $a_p(E):=(1+p)-\#\tilde{E}(\F_p)$ and $\tilde{E}$ is the reduced curve of $E$ modulo $p$. 

For $\ast\in\{\al, \be\}$, we write $\LL_{p, \ast}(t) \in \Q_p(\al)[[t]]$ for the cyclotomic $\ast$- $p$-adic $L$-function of $E$, which is defined with the image of Kato's zeta element for $E$ under Perrin-Riou's logarithm map. The definition of $\LL_{p, \ast}(t)$ is recalled in Subsection \ref{secpL}.

In the following, we assume $\ra(E/\Q)=1$. We introduce a $p$-adic analytic quantity
\[
\Sub_{\al,\be}:=\left\{\left(1-\frac{1}{\alpha}\right)^{-2}\LL_{p,\alpha}^{\prime} (0) - \left(1-\frac{1}{\beta}\right)^{-2}\LL_{p,\beta}^{\prime} (0)\right\}[\om, \varphi(\om)] \in \Q_p(\al).
\]
The factor $[\omega, \varphi(\omega)]$ $(\in \Z_p)$ is explained in Section \ref{secpL} later.
\begin{mainthm}[Theorem \ref{main1}]\label{intromain1}
Assume the conditions $(1)\sim (6)$ in Theorem \ref{main1}. Then
\[
v_p\left(L^{\prime}(E, 1)_{\alg} \cdot\Sub_{\al,\be}\right) \geq 1 \hspace{2mm}\Rightarrow\hspace{2mm} r(E[p]) \neq 0.
\] 
\end{mainthm}
\begin{rem}
The $p$-adic quantity $\Sub_{\al,\be}$ is inspired by the work of Kurihara--Pollack \cite{KuPo}, where they construct a rational point on an elliptic curve  numerically from a $p$-adic quantity like $\Sub_{\al,\be}$ when the elliptic curve has good supersingular reduction.
\end{rem}
\begin{rem}\label{introrefine}
We assume that $[\omega, \varphi(\omega)] \neq 0$ in the condition $(6)$ in Theorem \ref{main1}. In addition, as we will explain later, we have $v_p(\Sub_{\al,\be})\geq 0$ under the assumptions of Theorem \ref{intromain1}, and Theorem \ref{intromain1} refines (\ref{intromot}) by the $p$-adic quantity $\Sub_{\al,\be}$.
\end{rem}

Theorem \ref{intromain1} can detect the non-vanishing of $r(E[p])$ in Example \ref{introcounterexam}. 
\begin{example}\label{introexam}
We set $p=13$. Let $E$ be the elliptic curve defined by the equation $y^2+y=x^3+x^2$ as in Example \ref{introcounterexam}. The elliptic curve $E$ has good ordinary reduction at $13$. The conditions $(1), (3), (4), (5), (6)$ in Theorem \ref{main1} can be checked by numerical computations. The condition $(2)$ can be checked by the result of Kato \cite[Theorem 12.5]{Kato04} and Skinner--Urban \cite[Theorem 3.29]{SU14}. 

As we noted in Example \ref{introcounterexam}, we have $v_{13}(L^{\prime}(E, 1)_{\alg})=0$ in this case. In the proof of Theorem \ref{intromain1}, we show the equality of $p$-adic valuations
\[
v_p\left(\Sub_{\al,\be}\right) = v_p(\log_{\omega}(P)^2)+v_p(\al-\be)-2.
\]
Here $P$ is a generator of the torsion-free quotient $E(\Q)_{/\mathrm{tors}}$. The map $\log_{\om} : E(\Q_p) \to \Q_p$ is the formal logarithm associated with the N\'eron differential $\om$ for the fixed minimal model of $E/\Q_p$. We note that the valuation $v_{13}(\#\Sha(E/\Q))$ of the order of the Tate--Shafarevich group $\#\Sha(E/\Q)$ of $E$ is trivial since the $13$-part of the BSD conjecture is valid now as we noted in Example \ref{introcounterexam}.

We take $P$ as the point $(0,0) \in E(\Q)$. A computation shows that $v_{13}(\log_{\om}(P))= 3$, which implies $v_{13}\left(\Sub_{\al,\be} \right) = 4$ by the above equality. By Theorem \ref{intromain1}, we obtain the non-vanishing of $r(E[13])$.
\end{example}

\subsubsection{Result on division fields over imaginary quadratic fields}\label{introimag}
Next, we consider the case where the base field $F$ is an imaginary quadratic field $K$. Everywhere unramified $K$-rational points on $E$ also contribute the non-vanishing $r_K(E[p])\neq 0$. We can interpret this contribution in terms of the anticyclotomic $p$-adic $L$-function of $E$ over $K$ in this case.

Let $E$ be an elliptic curve over $\Q$ again, and we fix a minimal Weierstrass model of $E$ over $\Q$. In the following, we write $L(E/K, s)$ for the Hasse--Weil $L$-function of $E$ over $K$ and $\ra(E/K)$ for $\mathrm{ord}_{s=1}L(E/K, s)$. We consider the algebraic part 
\[
L^{\ast}(E/K, 1)_{\alg} := \frac{L^{\ast}(E/K, 1)}{|\mathrm{Disc}(K)|^{-1/2}\Omega_{E/K} \mathrm{Reg}(E/K)}
\]
of the leading coefficient $L^{\ast}(E/K, 1)$ of the Taylor expansion for $L(E/K, s)$ at $s=1$. Here $\mathrm{Disc}(K)$ is the discriminant of $K$, $\Omega_{E/K} := \int_{E(\C)}\om\wedge i\bar{\om}$ is a period of $E$ over $K$, and $\mathrm{Reg}(E/K)$ is the regulator of $E$ over $K$.

We assume that $K$ satisfies the following three conditions : 
\begin{itemize}
\item[(a)] $\mathrm{Disc}(K)$ is coprime to $pN$, where $N$ is the conductor of $E$.
\item[(b)] Every prime number which divides $N$ splits in $K$.
\item[(c)] The fixed prime number $p$ splits in $K$ as $p=\mathfrak{p}\bar{\mathfrak{p}}$. Here $\mathfrak{p}$ is the prime in $K$ above $p$ which is induced by the fixed embedding $\bar{\Q} \hookrightarrow \bar{\Q}_p$.
\end{itemize}

In this situation, we consider the anticyclotomic $p$-adic $L$-function $\LL_{\mathfrak{p}}^{\mathrm{BDP}}(t) \in W(\bar{\F}_p)[[t]]$ of $E$ over $K$, which is constructed by Bertolini--Darmon--Prasanna \cite{BDP13}, Brako\v{c}evi\'c \cite{Bra11} and Castella--Hsieh \cite{CH18}. Here $W(\bar{\F}_p)$ is the ring of Witt vectors over $\bar{\F}_p$. The square $L_{\mathfrak{p}}^{\mathrm{BDP}}(t):=(\LL_{\mathfrak{p}}^{\mathrm{BDP}}(t))^2$ interpolates central critical values of the Hasse--Weil $L$-function of $E/K$ twisted by certain infinite order characters. We assume that $\ra(E/K)=1$ in the following result.
\begin{mainthm}[Theorem \ref{main2}]\label{intromain2}
We assume that the conditions $(1)\sim (4)$ in Theorem \ref{main2} hold. Then we have 
\[
v_p\left(L^{\prime}(E/K, 1)_{\alg} \cdot L^{\mathrm{BDP}}_{\mathfrak{p}}(0)\right) \geq 1 \hspace{2mm}\Rightarrow \hspace{2mm} r_K(E[p]) \neq 0.
\]
\end{mainthm}
\begin{rem}
Since $\ra(E/K)=1$, we have $L^{\prime}(E/K, 1)_{\alg} \in \Q$ in Theorem \ref{intromain2} by the work of Gross--Zagier.
\end{rem}

We can observe many examples in which we can deduce $r_F(E[p])\neq 0$ only after considering the $p$-adic analytic quantities which we consider in Theorems \ref{intromain1} and \ref{intromain2}. Therefore, the author thinks it is important to consider such $p$-adic quantities besides complex $L$-values of $E$ in studying ideal class groups of division fields of elliptic curves. 

\subsection{Organization of this paper}
In Section \ref{unram prop}, we introduce some algebraic criteria for the non-vanishing of $r_F(E[p])$, which give some generalizations of the former results in \cite{PS21} and \cite{Da23}. In Section \ref{pL}, we first recall some basic notions on $p$-adic $L$-functions of elliptic curves. Then we introduce some known formulas on special values of $p$-adic $L$-functions, which we use in the proofs of Theorems \ref{intromain1} and \ref{intromain2}. In Section \ref{proof}, after stating Theorems \ref{intromain1} and \ref{intromain2} precisely, we conclude by proving them.

\section{Consideration of unramified cohomology groups}\label{unram prop}

\subsection{Basic notation}\label{notation}
First, we fix some notation. Let $F$ be a number field or a local field. We denote the absolute Galois group $\Gal(\overline{F}/F)$ of $F$ by $G_F$ for a fixed algebraic closure $\overline{F}$ of $F$.

For a number field $F$ and a place $v$ of $F$, $F_v$ denotes the completion of $F$ at $v$. We write $F^{\ur}_v$ for the maximal unramified abelian extension of $F_v$ in a fixed algebraic closure $\bar{F}_v$.

For a number field $F$ and a $G_F$-module $N$, we define the unramified subgroup of the Galois cohomology group $H^i(F, N):=H^i(G_F, V)$ as
\[
H_{\ur}^i(F, N):=\Ker\left(H^i(F, N) \xrightarrow{\prod_v\loc_v^{\ur}} \prod_{v\text{ : place}} H^i(F^{\ur}_v, N) \right).
\]
The map $\loc_v^{\ur}$ is the restriction of cohomology classes in $H^i(F, N)$ to the inertia subgroup $I_v$ in $G_F$ at a place $v$ of $F$.

Let $E$ be an elliptic curve over a local field $F$. We write $\mathcal{E}$ for the N\'eron model of $E$ over the ring of integers $\mathcal{O}_F$ of $F$. Let $E_0(F)$ be the subgroup of $E(F)=\mathcal{E}(\mathcal{O}_F)$ whose reduction lands in the identity component in the special fiber of $\mathcal{E}$. We also put $E_1(F)$ as the subgroup of $E(F)$ which has trivial reduction. The Tamagawa factor $c(E/F)$ of $E$ over $F$ is defined as $c(E/F):=[E(F) : E_0(F)]$.

Let $E$ be an elliptic curve over a number field $F$. We define the Tamagawa product of $E$ over $F$ as 
\[
\mathrm{Tam}(E/F):=\prod_{v\text{ : place}} c(E/F_v).
\]

For a prime number $p$, we write $F_E$ for the $p$-division field $F(E[p])$ of $E$ over $F$ as in Section \ref{intro}. We put the ideal class group of $F_E$ as $\cl(F_E)$ and $\cl(F_E)/p\cl(F_E)$ as $H_F$. We consider the semi-simplification
\begin{align}\label{semi}
H_F^{\semi} = \bigoplus_{M\hspace{1mm} :\hspace{1mm} \text{irr}} M^{\oplus r_F(M)}
\end{align}
of $H_F$ as an $\F_p[\Gal(F_E/F)]$-module. Here, for an irreducible $\F_p[\Gal(F_E/F)]$-module $M$, the non-negative integer $r_F(M)$ is the multiplicity of the $M$-component in $H_F^{\semi}$.
\subsection{Unramified cohomology groups and ideal class groups}
Let $F$ be a number field, $E$ an elliptic curve defined over $\Q$, and $p$ an odd prime number.

In the following, we study the multiplicity $r_{F}(E[p])$ of the $E[p]$-component in (\ref{semi}) when $E[p]$ is an irreducible $\Gal(F_E/F)$-module. For this purpose, it is important to investigate the unramified cohomology subgroup $H^1_{\ur}(F, E[p])$ of $H^1(F, E[p])$. 
\begin{prop}\label{unramcl}
Assume that $E[p]$ is an irreducible $\Gal(F_E/F)$-module. Then 
\[
H^1_{\ur}(F, E[p])\neq 0 \hspace{2mm}\Rightarrow\hspace{2mm} r_F(E[p])\neq 0.
\]
\end{prop}
\proofname\ By \cite[Lemma 2.2]{PS21}, the restriction map
\[
H^1(F, E[p]) \to \Hom_{\Gal(F_E/F)}(G_{F_E}, E[p])
\]
is injective. This map induces an injective map
\[
H^1_{\ur}(F, E[p]) \hookrightarrow \Hom_{\Gal(F_E/F)}(\Gal(F_E^{\ur}/F_E), E[p])\simeq \Hom_{\Gal(F_E/F)}(\cl(F_E), E[p]).
\]
Here we denote the maximal unramified abelian extension of $F_E$ by $F_E^{\ur}$ and use the isomorphism $\Gal(F_E^{\ur}/F_E) \simeq \cl(F_E)$ in class field theory. By the above injection, if we have $H^1_{\ur}(F, E[p]) \neq 0$, then there exists a non-trivial $\Gal(F_E/F)$-equivariant surjection $\cl(F_E)\twoheadrightarrow E[p]$ since $E[p]$ is irreducible. This implies the proposition. \hfill\qed

In the following, we introduce two results with which we can deduce the non-vanishing of $H^1_{\ur}(F, E[p])$, and the non-vanishing of $r_F(E[p])$ by Proposition \ref{unramcl}.

\subsection{Selmer groups and unramified cohomology groups}
\begin{prop}[]\label{PS result}
We assume the following conditions: 
\begin{itemize}
\item[(1)] The prime number $p$ does not divide $c(E/F_v)$ for any $v\nmid p$.
\item[(2)] For every place $v\mid p$ of $F$, $E(F_v)[p]=0$.
\end{itemize}
Then we have an inequality
\[
\dim_{\F_p} (H^1_{\ur}(F, E[p])) \geq \dim_{\F_p}(\Sel(F, E[p]))-[F:\Q].
\]
\end{prop}
\begin{rem}
When $F=\Q$, Proposition \ref{PS result} recovers part of \cite[Theorem 3.1]{PS21}.
\end{rem}
 
\proofname\ \ By \cite[Lemma 3.4]{AS02} and the assumption $(1)$, we can see that the localization
\[
\Sel(F, E[p]) \rightarrow H^1(F^{\ur}_v, E[p])
\]
of the $p$-Selmer group is the $0$-map for every place $v \nmid p$ of $F$. Then for the localization at $p$-adic places
\[
\prod_{v\mid p} \loc^{\ur}_v : \Sel(F, E[p]) \rightarrow \prod_{v\mid p} H^1(F^{\ur}_v, E[p]),
\]
we have
\begin{align}\label{SelinUr}
\Ker (\prod_{v\mid p} \loc^{\ur}_v)\subset H^1_{\ur}(F, E[p]).
\end{align}
We consider the localization $\loc_v : \Sel(F, E[p])\to H^1(F_v, E[p])$ at $v\mid p$. By the definition of the Selmer group and the assumption $(2)$, we have $\mathrm{Im}(\loc_v) \hookrightarrow E(F_v)/pE(F_v) \simeq \F_p^{[F_v : \Q_p]}$. We obtain
\[
\dim_{\F_p}\mathrm{Im} (\prod_{v\mid p} \loc^{\ur}_v ) \leq \dim_{\F_p}\mathrm{Im} (\prod_{v\mid p} \loc_v ) \leq \sum_{v\mid p} [F_v : \Q_p]=[F:\Q],
\]
and
\[
\dim_{\F_p}\Ker (\prod_{v\mid p} \loc^{\ur}_v) \geq \dim_{\F_p}(\Sel(F,E[p]))-[F:\Q].
\]
The last inequality and (\ref{SelinUr}) imply the proposition.\hfill\qed

As an immediate consequence of Propositions \ref{unramcl} and \ref{PS result}, we have the following.
\begin{cor}
We assume the conditions in Propositions \ref{unramcl} and \ref{PS result}. Then
\[
\dim_{\F_p}(\Sel(F, E[p])) \geq [F:\Q]+1 \hspace{2mm}\Rightarrow \hspace{2mm}r_F(E[p]) \neq 0.
\]
\end{cor}

\subsection{Rational points and unramified cohomology groups}
Next, we give a condition for rational points of elliptic curves to yield unramified cohomology classes via the Kummer map. 
\begin{prop}[]\label{my result}
We fix a place $v\mid p$ of $F$. We assume the following conditions: 
\begin{itemize}
\item[(1)] The prime number $p$ does not divide $c(E/F_v)$ for any $v\nmid p$.
\item[(2)] The extension $F/\Q$ is a Galois extension.
\item[(3)] The ramification index $e(F_v/\Q_p)$ of $F_v/\Q_p$ satisfies $\displaystyle e(F_v/\Q_p)<p-1$. 
\end{itemize}
We further assume that $E$ is defined by a minimal model over $F_v$. If there exists a rational point $Q \in (E(F)\cap E_1(F_v))\setminus pE(F)$ such that 
\[
v_p(\log_{\omega}(Q)) \geq \frac{e(F_v/\Q_p)+1}{e(F_v/\Q_p)}, 
\]
then we have $H^1_{\ur}(F, E[p])\neq 0$. Here $\log_{\om} : E(F_v) \to F_v$ is the formal group logarithm associated to the N\'eron differential $\om$ for the fixed minimal model of $E$ over $F_v$.
\end{prop}

\begin{rem}
When $F=\Q$, Proposition \ref{my result} recovers \cite[Theorem 4.5]{Da23}.
\end{rem}

\proofname \ \ We consider the image of the rational point $Q$ under the Kummer map 
\[
\kappa : E(F)/pE(F) \hookrightarrow \Sel(F, E[p]).
\]
 
We show that $\kappa(Q)$ gives an element in $H^1_{\ur}(F, E[p])$.

By \cite[Lemma 3.4]{AS02} and the assumption $(1)$, we can see that elements in $\Sel(F, E[p])$ are trivial when restricted to the inertia subgroup at any place $v^{\prime}\nmid p$ of $F$. We show that $\kappa(Q)$ is trivial when restricted to the inertia subgroup at any place $v^{\prime} \mid p$ of $F$. By the assumption $(2)$, it suffices to show the triviality for the fixed place $v\mid p$.

Let $\mathcal{O}_{F_v}$ be the ring of integers of $F_v$ and $\pi$ a uniformizer of $\mathcal{O}_{F_v}$. The condition in Proposition \ref{my result} on the $p$-adic valuation of $\log_{\om}(Q)$ implies that $\log_{\om}(Q)$ is a multiple of $p$ in the maximal ideal $\pi \mathcal{O}_{F_v}$. This further implies that $Q \in pE_1(F_v)$ since $\log_{\om}$ gives an isomorphism $\log_{\om} : E_1(F_v) \xrightarrow{\simeq} \pi \mathcal{O}_{F_v}$ by the assumption $(3)$. See \cite[Theorem 6.4 (b)]{Sil} for this isomorphism. We consider the diagram
\[
  \xymatrix{
      E(F)/pE(F) \ar[r]^{\kappa}\ar[d] & \Sel(F, E[p])\ar[d]^{\loc_v}\\
  E(F_v)/pE(F_v)\ar[r]^{\kappa_v} &  H^1(F_v,  E[p]),
  }
\]
where $\kappa_v$ is the local Kummer map and the left vertical arrow is induced by the inclusion. This diagram implies $\loc_v(\kappa(Q))=0$ and $\kappa(Q) \in H^1_{\ur}(F, E[p])\setminus\{0\}$.\hfill\qed

As an immediate consequence of Propositions \ref{unramcl} and \ref{my result}, we have the following. 
\begin{cor}
We assume the conditions in Propositions \ref{unramcl} and \ref{my result}. For a fixed place $v\mid p$ of $F$, if there exists a rational point $Q \in (E(F)\cap E_1(F_v))\setminus pE(F)$ which satisfies 
\[
v_p(\log_{\omega}(Q)) \geq \frac{e(F_v/\Q_p)+1}{e(F_v/\Q_p)}, 
\]
then we have $r_F(E[p]) \neq 0$.
\end{cor}

\section{$p$-adic $L$-functions of elliptic curves and special value formulas}\label{pL}
In this section, we review $p$-adic $L$-functions of elliptic curves and introduce some known special value formulas of them for later use.

\subsection{Cyclotomic $p$-adic $L$-functions of elliptic curves}\label{secpL}
We fix an odd prime number $p$. Let $E$ be an elliptic curve over $\Q$ which has good reduction at $p$. We fix a minimal Weierstrass model of $E/\Q$. We first introduce a setup for the crystalline Dieudonné module $\D(V)$ of $V:=T \tens{\Z_p}{\Q_p}$, where $T:=T_p(E)$ denotes the $p$-adic Tate module of $E$.

Let $\al, \beta \in \bar{\Q}_p$ be the two roots of the Hecke polynomial $X^2-a_p(E)X+p$ which satisfies $v_p(\al)\leq v_p(\be)$, where $a_p(E):=(p+1)-\#\tilde{E}(\F_p)$. We set $L:=\Q_p(\alpha)$. 

We write $\varphi$ for the Frobenius endomorphism on $\D(V)$, whose eigenvalues are $\al^{-1}, \be^{-1} \in L$. For $\ast\in\{\al,\be\}$, let $D_{\ast}$ be the eigenspace of $\varphi$ in the $L$-vector space $\D(V)_L:=\D(V)\tens{\Q_p}{L}$ associated with the eigenvalue $\ast^{-1}$.

Let $\om \in H^1_{\mathrm{dR}}(E/\Q_p)$ be the N\'eron differential of $E$ associated with the fixed minimal model of $E/\Q_p$, which we can think as an element in $\D(V)$ via the comparison isomorphism in the $p$-adic Hodge theory.

We write $[\cdot,\cdot] : \D(V) \times\D(V) \to \Q_p$ for the alternating pairing induced by the cup product on $H^1_{\mathrm{dR}}(E/\Q_p)$, which is normalized as $[\om, \eta]=1$. Here $\eta:=x\om \in H^1_{\mathrm{dR}}(E/\Q_p)$ and $x$ denotes the $x$-coordinate of the defining polynomial of $E$. We extend $[\cdot,\cdot]$ on $\D(V)_L$ linearly and denote by the same notation.

From now on, we assume that $\{\om, \varphi(\om)\}$ $\subset \D(V)$ are {\it linearly independent}. Let $\om_{\al}$ and $\om_{\be}$ be the projections of $\omega \in \D(V)_L$ on $D_{\al}$ and $D_{\be}$ respectively:
\[
\om_{\al}:= \frac{\al\be}{\al-\be}(\be^{-1}\omega-\varphi(\om)) \in D_{\al},\hspace{3mm} \om_{\be}:=- \frac{\al\be}{\al-\be}(\al^{-1}\omega-\varphi(\om))\in D_{\be}.
\]
\begin{rem}
When $E$ has good supersingular reduction at $p$, it is known that $[\om, \varphi(\om)]\in \Z_p^{\times}$. In fact, $\{\om, \varphi(\om)\}$ are linearly dependent if and only if $E$ has good ordinary reduction at $p$ and has complex multiplication. In such a case, in other words, $\om$ is an eigenvector for $\varphi$, then one of $\om_{\al}$ and $\om_{\be}$ is $0$ and the other is $\om$. 
\end{rem} 
  
Next, we recall some tools in Perrin-Riou's theory developed in \cite{PR93} and \cite{PR94}. Let $\Q_{\infty}/\Q$ be the cyclotomic $\Z_p$-extension of $\Q$ and $\Q_{\infty, p}$ be the completion of $\Q_{\infty}$ at the unique prime ideal above $p$. We put the Galois group of $\Q_{\infty}/\Q$ as $\Gamma$, which can be identified with the Galois group of $\Q_{\infty,p}/\Q_p$ via the restriction.

Put 
\[
\hh :=\left\{f(t):=\sum_{n=0}^{\infty} a_nt^n \in\Q_p[[t]] \middle| \exists h\in \Z_{>0}\hspace{2mm} \text{s.t} \lim_{n\to\infty} \frac{|a_n|_p}{n^h}=0\right\},
\]
where $|\ast|_p$ is the $p$-adic absolute value defined as $|x|_p=p^{-v_p(x)}$. In the following, we fix a topological generator $\gamma$ of $\Gamma$ and put
\[
\hh(\Gamma) := \left\{f(\gamma-1) \middle| f(t) \in \hh\right\}.
\]
 
Let $\HH^1(\Q_{\infty,p}, T):=\varprojlim H^1(\Q_{n,p}, T)$ be the local Iwasawa cohomology group for $T$. Here $\Q_{n, p}$ is the completion of the $n$-th layer $\Q_n$ of $\Q_{\infty}/\Q$ at the unique $p$-adic place. Let
\[
\LL_V : \HH^1(\Q_{\infty,p}, T) \to \hh(\Gamma)\tens{\Q_p}{\D(V)}
\] 
be Perrin-Riou's big logarithm map, which is the same as $-(\Omega^{\varepsilon}_{V,0})^{-1}$ in \cite{PR93} for a fixed compatible system $\varepsilon = (1, \zeta_p, \zeta_p^2, \ldots , )$ of root of unities of $p$-power order. For $\bm{x} \in \HH^1(\Q_{\infty,p}, T)$, we denote its image under the map $\LL_V$ by $\LL_{V,\bm{x}}$.

Hereafter, we assume that $E[p]$ is irreducible as a $\Gal(\Q_E/\Q)$-module, where $\Q_E:=\Q(E[p])$. Let 
\[
\z^{\mathrm{Kato}}=(z_n)_n \in \HH^1(\Q_{\infty}, T)
\]
be Kato's zeta element introduced in \cite{Kato04}. We denote the $n$-th layer of $\z^{\mathrm{Kato}}$ by $z_n \in H^1(\Q_n, T)$. We often denote the localization of $\z^{\mathrm{Kato}}$ and $z_n$ at $p$ by the same notation if no confusion occurs. The element $z_n$ has the following properties (Kato's explicit reciprocity law): 
\begin{itemize}
\item[(1)] For any $n>0$ and faithful character $\psi : \Gal(\Q_n/\Q) \to \bar{\Q}_p^{\times}$, we have 
\[
\sum_{\sigma \in \Gal(\Q_n/\Q)}\psi(\sigma)\exp_V^{\ast}(\sigma(z_n)) = \frac{L(E, \psi, 1)}{\Omega_E}\omega.
\]
Here $L(E, \psi, s)$ is the Hasse--Weil $L$-function of $E$ twisted by $\psi$ and $\exp_V^{\ast} : H^1(\Q_{n, p}, V) \to \Q_{n, p}\tens{\Q_p}{\mathrm{Fil}^0\D(V)}$ is the dual exponential map.
\item[(2)] 
\[
\exp^{\ast}_V(z_0) = (1-a_p(E)p^{-1}+p^{-1})\frac{L(E, 1)}{\Omega_E}\omega.
\]
\end{itemize}

Let $\LL_{p}(t)$ be the $\D(V)$-valued power series which corresponds to the element $\LL_{V, \z^{\mathrm{Kato}}} \in \hh(\Gamma)\tens{\Q_p}{\D(V)}$ under the isomorphism $\hh(\Gamma)\tens{\Q_p}{\D(V)} \simeq \hh\tens{\Q_p}{\D(V)}$. We define $\LL_{p, \ast}(t) \in L[[t]]$ for $\ast \in \{\al, \be\}$, using the decomposition 
\[
\hh\tens{\Q_p}{\D(V)_L} = \hh\tens{\Q_p}{(L\omega_{\al}\oplus L\omega_{\be})} \simeq (\hh\tens{\Q_p}L)^{\oplus 2}\hspace{2mm}:\hspace{2mm}\LL_{p}(t) \mapsto (\LL_{p, \al}(t), \LL_{p, \be}(t)).
\]
The power series $\LL_{p, \ast}(t)$ is the so-called $\ast$-\hspace{0.5mm} $p$-adic $L$-function of $E$. These functions interpolate special values of the $L$-function of $E/\Q$ twisted by characters of $\Gal(\Q_n/\Q)$. For example, the $\al$-\hspace{0.5mm}$p$-adic $L$-function $\LL_{p, \al}(t)$ satisfies
\[
\LL_{p, \al}(0) = 
\left(1-\dfrac{1}{\al}\right)^2\dfrac{L(E, 1)}{\Omega_E},
\]
which is a consequence of the descent property \cite[Proposition 2.1.4]{PR93} of Perrin-Riou's map $\LL_V$ and Kato's explicit reciprocity law.

\subsection{Formulas concerning cyclotomic $p$-adic $L$-functions} 
A $p$-adic analogue of BSD conjecture for the $\al$-\hspace{0.5mm}$p$-adic $L$-function $\LL_{p, \al}(t)$ of $E$ is originally formulated by Mazur--Tate--Teitelbaum and stated in \cite{Co04}, including the good supersingular case. We state it only for the case where $\ra(E/\Q)=1$.

\begin{conj}[{a special case of \cite[Conjecture 0.12 (i)]{Co04}}]\label{paBSD}
Suppose $\ra(E/\Q)=1$. 
\[
\LL_{p, \al}^{\prime}(0) = \frac{1}{\log_p(\chi_{\cyc}(\gamma))}\left(1-\dfrac{1}{\al}\right)^2\dfrac{\#\Sha(E/\Q)\cdot \mathrm{Tam}(E/\Q) \cdot \langle P, P \rangle_{p,\al,L}}{(\#E(\Q)_{\mathrm{tors}})^2}.
\]
Here $\chi_{\cyc} : \Gamma \to 1+p\Z_p$ is the $p$-adic cyclotomic character and
\[
\langle \cdot, \cdot \rangle_{p,\al, L}:=\langle \cdot, \cdot \rangle_{p,L} : (E(\Q)\tens{\Z}{L}) \times (E(\Q)\tens{\Z}{L}) \to L
\]
is the $p$-adic height pairing associated to the unit root $\al$. The rational point $P \in E(\Q)$ is a generator of the torsion-free quotient $E(\Q)_{/\mathrm{tors}}$.
\end{conj}

It is known that the Iwasawa main conjecture (IMC) for $E$ and $p$ implies the $p$-part of Conjecture \ref{paBSD} under certain conditions. 
\begin{prop}[a special case of {\cite[Theorem $2^{\prime}$]{Sch85} and \cite[Proposition 3.4.6]{PR93}}]\label{IMCtopaBSD}\mbox{}

Suppose $\ra(E/\Q)=1$. We further assume the following conditions:
\begin{itemize}
\item[(1)] If $p\nmid a_p(E)$, then the IMC in \cite{Maz72, Gr99} holds and the $p$-adic height pairing $\langle \cdot, \cdot \rangle_{p,L}$ associated to $\al$ is non-trivial.
\item[(2)] If $p\mid a_p(E)$, then the IMC in \cite[Conjecture CP(V), Section 3.4]{PR93} holds. 
\end{itemize}
Then the $p$-part of  Conjecture \ref{paBSD} holds:  
\[
v_p\left(\LL_{p, \al}^{\prime}(0)\right) = v_p\left(\frac{1}{\log_p(\chi_{\cyc}(\gamma))}\left(1-\dfrac{1}{\al}\right)^2\dfrac{\#\Sha(E/\Q)\cdot \mathrm{Tam}(E/\Q) \cdot \langle P, P \rangle_{p, L}}{(\#E(\Q)_{\mathrm{tors}})^2}\right).
\] 
\end{prop}
\begin{rem}
The IMC in \cite{Maz72, Gr99} is proved by Kato and Skinner--Urban under certain mild assumptions. The IMC \cite[Conjecture CP(V), Section 3.4]{PR93} is equivalent to the conjecture \cite[Conjecture, Section 1]{Kob03} formulated by Kobayashi. Recently, the proof of \cite[Conjecture, Section 1]{Kob03} was announced in \cite{BSW24} when $a_p(E)=0$ and $E$ has a square-free conductor.
\end{rem}

Here we introduce a formula for a specific value of the $p$-adic height pairing $\langle \cdot, \cdot \rangle_{p, L}$. Recall that we write $z_0$ for the $0$-th layer of $\z^{\mathrm{Kato}}$. When $\ra(E/\Q)=1$, we have 
$z_0 \in \Sel(\Q, T)\simeq E(\Q)\tens{\Z}{\Z_p}$ by Kato's explicit reciprocity law.

\begin{prop}\label{gen.rubin}
We assume $\ra(E/\Q)=1$. For any $x \in E(\Q)\tens{\Z}{L}$, we have 
\[
\langle z_0, x \rangle_{p, L} = \left(1-\frac{1}{\al}\right)^{-1}\left(1-\frac{1}{\beta}\right) \log_{\om}(x)\cdot \LL_{p,\al}^{\prime}(0)\cdot\log_p(\chi_{\cyc}(\gamma)).
\]
\end{prop}
\begin{rem}
When $E$ has good ordinary reduction at $p$, this formula is shown by Rubin. This is just a consequence of Perrin--Riou's results \cite[Proposition 2.2.4, Proposition 2.3.4 (bis)]{PR93}, but we explain the proof briefly below.
\end{rem}

\proofname \ First, since we assume $\ra(E/\Q)=1$, we have $\LL_p(0)=0$ by \cite[Proposition 2.1.4]{PR93} and Kato's explicit reciprocity law. We take $D_{\be}=L\cdot \om_{\be}$ as the complement subspace of $\mathrm{Fil}^0\D(V)$ in $\D(V)$. The projection of $(1-p^{-1}\varphi^{-1})(1-\varphi)^{-1} \LL_p^{\prime}(0)$ onto $\mathrm{Fil}^0\D(V)$ with respect to $D_{\be}$ is calculated in \cite[Proposition 2.2.4]{PR93} and denoted by $\delta_V(\z^{\mathrm{Kato}}) \in \mathrm{Fil}^0\D(V)$. Then we have
\[
\log_p(\chi_{\cyc}(\gamma))[(1-p^{-1}\varphi^{-1})(1-\varphi)^{-1} \LL_p^{\prime}(0), \omega_{\beta}] = [\delta_V(\z^{\mathrm{Kato}}), \om_{\be}].
\] 
For the left-hand side, we have 
\[
[(1-p^{-1}\varphi^{-1})(1-\varphi)^{-1} \LL_p^{\prime}(0), \omega_{\beta}] = \left(1-\frac{1}{\be}\right)\left(1-\frac{1}{\al}\right)^{-1}\LL_{p, \al}^{\prime}(0)\cdot [\om_{\al}, \om_{\be}].
\]
On the other hand, \cite[Proposition 2.3.4 (bis)]{PR93} gives a relation between values of the $p$-adic height pairing and the cup product pairing as 
\[
\langle z_0, x \rangle_{p, L} = [\delta_V(\z^{\mathrm{Kato}}), \log_V(x)].
\]
Here, $\log_V : E(\Q_p) \tens{\Z}{L} \to \D(V)_L$ is the Bloch--Kato logarithm map. We can show that
\[
\log_{\omega}(x)[\delta_V(\z^{\mathrm{Kato}}), \om_{\be}] = [\delta_V(\z^{\mathrm{Kato}}), \log_V(x)]\cdot [\om_{\al}, \om_{\be}]
\]
by the relation $[\om, \log_V(x)]=\log_{\omega}(x)$ and the definitions of $\om_{\al}, \om_{\be}$. Thus the proposition follows. \hfill\qed

\subsection{Anticyclotomic $p$-adic $L$-functions of elliptic curves}\label{antipL}
Let $E$ be an elliptic curve $E/\Q$ which has good reduction at an odd prime number $p$. We take $K$ as an imaginary quadratic field which satisfies the following conditions.
\begin{itemize}
\item[(a)] The discriminant of $K$ is coprime to $pN$, where $N$ is the conductor of $E$.
\item[(b)] Every prime number which divides $N$ splits in $K$.
\item[(c)] The fixed prime number $p$ splits in $K$ as $p=\mathfrak{p}\bar{\mathfrak{p}}$.
\end{itemize}
In the condition (c), we write $\mathfrak{p}$ for the $p$-adic place of $K$ which is induced by the fixed embedding $\bar{\Q} \hookrightarrow \bar{\Q}_p$.

We denote the anticyclotomic $\Z_p$-extension of $K$ by $K_{\infty}$ and its Galois group over $K$ by $\Gamma_K$. As in the cyclotomic case, we fix a topological generator $\gamma_K$ of $\Gamma_K$. 

Let $\LL_{\mathfrak{p}}^{\mathrm{BDP}} \in W(\bar{\F}_p)[[\Gamma_K]]$ be the anticyclotomic $p$-adic $L$-function constructed by Bertolini--Darmon--Prasanna \cite{BDP13}, Brako\v{c}evi\'c \cite{Bra11} and Castella--Hsieh \cite{CH18}. We can think of this function as a power series $\LL_{\mathfrak{p}}^{\mathrm{BDP}}(t)$ in $W(\bar{\F}_p)[[t]]$ by the isomorphism $W(\bar{\F}_p)[[\Gamma_K]] \simeq W(\bar{\F}_p)[[t]]$ which sends $\gamma_K \mapsto 1+t$. We put $L_{\mathfrak{p}}^{\mathrm{BDP}}(t):=(\LL_{\mathfrak{p}}^{\mathrm{BDP}}(t))^2$, which interpolates central critical values of the complex $L$-function of $E/K$ twisted by certain infinite order characters of $\Gamma_K$.

In \cite[Conjecture 1.1]{CHKLL23}, a $p$-adic analogue of the BSD conjecture for $L_{\mathfrak{p}}^{\mathrm{BDP}}(t)$ is formulated. We state the conjecture only the case where $\ra(E/K)=1$.
\begin{conj}[{a special case of \cite[Conjecture 1.1]{CHKLL23}}]\label{antipaBSD}
Assume that $\ra(E/K)=1$ and $E(K)[p]=0$. We have an equality 
\[
L_{\mathfrak{p}}^{\mathrm{BDP}}(0) = \left(1-a_p(E)p^{-1}+p^{-1}\right)^2 \cdot \#\Sha(E/K) \cdot \log_{\om}(P)^2 \cdot \mathrm{Tam}(E/\Q)^2,
\]
where $P$ is a generator of the torsion-free quotient $E(K)_{/\mathrm{tors}}$.
\end{conj}
\begin{rem}
Since $\ra(E/K)=1$, we have $\ral(E/K)=1$ and the finiteness of $\#\Sha(E/K)$ by the result of Gross--Zagier and Kolyvagin.
\end{rem}
In \cite{CHKLL23}, an analogous result of Proposition \ref{IMCtopaBSD} for this anticyclotomic setting is proved.
\begin{prop}[{a special case of \cite[Theorem 5.1]{CHKLL23}}]\label{antiIMCtopaBSD}
Suppose $\ra(E/K)=1$. We further assume that the anticyclotomic Iwasawa main conjecture (IMC) \cite[Conjecture 2.4]{CHKLL23} holds for $E$ and $p$. Then the $p$-part of the Conjecture \ref{antipaBSD} holds: 
\[
v_p\left(L_{\mathfrak{p}}^{\mathrm{BDP}}(0)\right) = v_p\left(\left(1-a_p(E)p^{-1}+p^{-1}\right)^2 \cdot \#\Sha(E/K) \cdot \log_{\om}(P)^2 \cdot  \mathrm{Tam}(E/\Q)^2\right).
\]
\end{prop}
\begin{rem}
The anticyclotomic IMC \cite[Conjecture 2.4]{CHKLL23} is proved in  
\cite[Theorem B]{BCK21} when $p\nmid a_p(E)$ and in \cite[Corollary 7.2]{CHKLL23} when $a_p(E)=0$ under certain mild conditions.
\end{rem}

\section{Main results and proofs}\label{proof}
In this section, we state Theorems \ref{intromain1} and \ref{intromain2} in Section \ref{introresults} precisely and prove them. 

\subsection{Result on division fields over $\Q$}
Let $E$ be an elliptic curve over $\Q$ and $p$ an odd prime number. We assume that $E$ has good reduction at $p$.

\begin{defn}\label{defS}
\[
\Sub_{\al,\be}:=\left\{\left(1-\frac{1}{\alpha}\right)^{-2}\LL_{p,\alpha}^{\prime} (0) - \left(1-\frac{1}{\beta}\right)^{-2}\LL_{p,\beta}^{\prime} (0)\right\}[\om, \varphi(\om)] \in L(=\Q_p(\al)).
\]
\end{defn}

\begin{thm}[Theorem \ref{intromain1} in Section \ref{introresults}]\label{main1}
We assume the following conditions: 
\begin{itemize}
\item[(1)] The elliptic curve $E$ has good reduction at $p$ and $\ra(E/\Q)=1$.
\item[(2)] When $p\nmid a_p(E)$ (resp. $p\mid a_p(E)$), the Iwasawa main conjecture in $(1)$ (resp. $(2)$) in Proposition \ref{IMCtopaBSD} holds for $E$ and $p$.
\item[(3)] If $p\nmid a_p(E)$, then the $p$-adic height pairing $\langle \cdot, \cdot \rangle_{p, \al, L}:=\langle \cdot, \cdot \rangle_{p, L}$ associated to the unit root $\alpha$ is non-trivial.  
\item[(4)] The prime number $p$ does not divide $\# E(\F_p)\cdot \mathrm{Tam}(E/\Q)$.
\item[(5)] As a $\Gal(\Q_E/\Q)$-module, $E[p]$ is irreducible.
\item[(6)] The elements $\om, \varphi(\om) \in \D(V)$ are linearly independent, in other words, $[\om, \varphi(\om)] \neq 0$.
\end{itemize}
Then we have 
\begin{align}\label{imp}
v_p\left(L^{\prime}(E,1)_{\alg} \cdot \Sub_{\al,\be}\right) \geq 1  \hspace{2mm}\Rightarrow\hspace{2mm} r(E[p]) \neq 0.
\end{align}
Here $L^{\prime}(E,1)_{\alg}$ is the algebraic part of $L^{\prime}(E,1)$ defined in Section \ref{intro}.  
\end{thm}
\begin{rem}\label{an alg}
By the work of Gross-Zagier and Kolyvagin, we have $\ral(E/\Q):=\mathrm{rank}_{\Z}(E(\Q))=1$ in Theorem \ref{main1} since $\ra(E/\Q)=1$.
\end{rem}

Now we prove Theorem \ref{main1}. We first introduce the following two lemmas.
\begin{lem}\label{subst} We assume the conditions $(1), (4)$ in Theorem \ref{main1}.
Put 
 \[
 \Sub_p  := [(1-p^{-1}\varphi^{-1})(1-\varphi)^{-1}\LL_p^{\prime}(0), \hspace{2mm}\omega]\hspace{1mm} \in \Q_p.
 \]
Then we have 
\[
v_p\left(\Sub_{\alpha,\beta}\right) = v_p\left(\Sub_p\right) + v_p\left(\al-\be\right).
\]
\end{lem}
 \proofname\ \ Recall that we take a decomposition $\om=\om_{\al}+\om_{\be}$ in Subsection \ref{secpL}, where $\om_{\al}, \om_{\be} \in \D(V)_L$ are eigenvectors of $\varphi$ for the eigenvalues $\al^{-1}$ and $\be^{-1}$ respectively. We have
\begin{align*}
\Sub_p & = [(1-p^{-1}\varphi^{-1})(1-\varphi)^{-1}\LL_p^{\prime}(0),\hspace{2mm}\omega]\\
& = [(1-p^{-1}\varphi^{-1})(1-\varphi)^{-1}(\LL_{p,\alpha}^{\prime}(0)\omega_{\alpha}+\LL_{p,\beta}^{\prime}(0)\omega_{\beta}),\hspace{2mm}\omega_{\alpha}+\omega_{\beta}]\\
& = \left[\left(1-\frac{1}{\beta}\right)\left(1-\frac{1}{\alpha}\right)^{-1}\LL_{p,\alpha}^{\prime}(0)\omega_{\alpha}+\left(1-\frac{1}{\alpha}\right)\left(1-\frac{1}{\beta}\right)^{-1}\LL_{p,\beta}^{\prime}(0)\omega_{\beta},\hspace{2mm}\omega_{\alpha}+\omega_{\beta}\right]\\
& = \left(\left(1-\frac{1}{\beta}\right)\left(1-\frac{1}{\alpha}\right)^{-1}\LL_{p,\alpha}^{\prime}(0)-\left(1-\frac{1}{\alpha}\right)\left(1-\frac{1}{\beta}\right)^{-1}\LL_{p,\be}^{\prime}(0)\right)\cdot[\omega_{\al}, \omega_{\be}].
\end{align*} 
The value $[\omega_{\al}, \omega_{\be}]$ can be computed as $
[\omega_{\al}, \omega_{\be}] = \dfrac{p}{\al-\be}[\omega, \varphi(\omega)]
$.
Then
\begin{align*}
\Sub_p & = \left(1-\frac{1}{\al}\right)\left(1-\frac{1}{\be}\right)\Sub_{\al,\be}\cdot \frac{p}{\al-\be}\\
& = \dfrac{1-a_p(E)+p}{\al-\be}\Sub_{\al,\be}.
\end{align*}
The factor $1-a_p(E)+p$ is a $p$-adic unit due to the assumption (4) in Theorem \ref{main1}, which implies the claim. \hfill\qed

\begin{lem}\label{KEY}
We assume the conditions $(1), (2), (3), (4), (6)$ in Theorem \ref{main1}. We have 
\[
v_p(L^{\prime}(E,1)_{\alg}\cdot \Sub_p) = v_p(\#\Sha(E/\Q)^2\cdot \log_{\om}(P)^2)-2,
\]
where $P \in E(\Q)$ is a generator of $E(\Q)_{/\mathrm{tors}}$.
\end{lem}

\proofname\ \ Since we assume $\ra(E/\Q)=1$ in the condition $(1)$, we have $z_0 \in \Sel(\Q, T) \simeq E(\Q) \tens{\Z}\Z_p$ by Kato's explicit reciprocity law, and also $\LL_p(0)=0$ by \cite[Proposition 2.1.4]{PR93}. By the formula \cite[Proposition 2.2.2]{PR93} for the derivative $\LL_p^{\prime}(0)$, we have
\[
(1-p^{-1}\varphi^{-1})(1-\varphi)^{-1}\LL_p^{\prime}(0) \equiv \log_p(\chi_{\cyc}(\gamma))^{-1}\cdot \log_V(z_0) \mod \mathrm{Fil}^0\D(V).
\]
From this equality,
\begin{align}\label{c of Pe}
\Sub_p& =[(1-p^{-1}\varphi^{-1})(1-\varphi)^{-1}\LL_p^{\prime}(0),\omega]\nonumber\\
& = [\log_p(\chi_{\cyc}(\gamma))^{-1}\cdot \log_V(z_0), \omega]\nonumber\\
& = -\log_p(\chi_{\cyc}(\gamma))^{-1}\log_{\omega} (z_0).
\end{align}

While, putting $z_0 \in E(\Q) \tens{\Z}\Z_p$ into $x$ in Proposition \ref{gen.rubin}, we have 
\begin{align*}
& \langle z_0, z_0 \rangle_{p,L} = \left(1-\dfrac{1}{\alpha}\right)^{-1}\left(1-\dfrac{1}{\beta}\right)\log_{\omega}(z_0)\cdot \LL_{p,\alpha}^{\prime}(0)\cdot\log_p(\chi_{\cyc}(\gamma))\\
\iff & \log_{\omega}(z_0) = \left(1-\dfrac{1}{\alpha}\right)^{-1}\left(1-\dfrac{1}{\beta}\right) \LL_{p,\alpha}^{\prime}(0)\cdot\log_p(\chi_{\cyc}(\gamma))\cdot \dfrac{\log_{\omega}(z_0)^2}{\langle z_0, z_0 \rangle_{p,L}}.
\end{align*}
Since $\ral(E/\Q)=1$ by Remark \ref{an alg}, we have
\begin{align}\label{c of Ru}
\log_{\omega}(z_0) = \left(1-\dfrac{1}{\alpha}\right)^{-1}\left(1-\dfrac{1}{\beta}\right) \LL_{p,\alpha}^{\prime}(0)\cdot\log_p(\chi_{\cyc}(\gamma))\cdot \dfrac{\log_{\omega}(P)^2}{\langle P, P \rangle_{p,L}}
\end{align}
for a generator $P$ of $E(\Q)_{/\mathrm{tors}}$. By the assumptions (2) and (3) in Theorem \ref{main1} and Proposition \ref{IMCtopaBSD}, the $p$-part of Conjecture \ref{paBSD} for $\LL_{p,\alpha}$ is now valid, which asserts
 \begin{align}\label{c of paBSD}
 &v_p\left(\left(1-\dfrac{1}{\alpha}\right)^{-1}\left(1-\dfrac{1}{\beta}\right) \LL_{p,\alpha}^{\prime}(0) \cdot\log_p(\chi_{\cyc}(\gamma))\right)\nonumber\\
  = & v_p\left((1-a_p(E)p^{-1}+p^{-1})\cdot\dfrac{\#\Sha(E/\Q)\cdot \mathrm{Tam}(E/\Q)}{\#E(\Q)_{\mathrm{tors}}^2} \langle P, P \rangle_{p,L}\right)\nonumber\\
  = & v_p\left(\#\Sha(E/\Q)\cdot \langle P, P \rangle_{p,L}\right)-1.
 \end{align}
Here we note that the $p$-adic valuations of the factors $\mathrm{Tam}(E/\Q)$, $\#E(\Q)_{\mathrm{tors}}^2$ and $1-a_p(E)p^{-1}+p^{-1}$ are $0$, $0$ and $-1$, respectively by the assumption (4) in Theorem \ref{main1}. Thus (\ref{c of Pe}), (\ref{c of Ru}) and (\ref{c of paBSD}) give an equality
\begin{align}\label{summary}
v_p(\Sub_p) =  v_p\left(\#\Sha(E/\Q)\log_{\omega}(P)^2\right)-2.
\end{align}

Due to \cite[Corollary 1.8]{PR87} when $p\nmid a_p(E)$ and \cite[Corollary 1.3 (ii)]{Kob13} when $p\mid a_p(E)$, the $p$-part of Conjecture \ref{paBSD} is now equivalent to the $p$-part of the BSD conjecture by the assumption $\ra(E/\Q)=1$ and $(3)$ in Theorem \ref{main1}. Then we have
\[
v_p(L^{\prime}(E,1)_{\alg})=v_p\left(\frac{\#\Sha(E/\Q)\cdot\mathrm{Tam}(E/\Q)}{\#E(\Q)^2_{\mathrm{tors}}}\right)= v_p(\#\Sha(E/\Q)).
\]
Combining this equality with (\ref{summary}), we obtain
\[
v_p(L^{\prime}(E,1)_{\alg}\cdot\Sub_p) = v_p\left(\#\Sha(E/\Q)^2\log_{\omega}(P)^2\right)-2.
\]
This completes the proof.\hfill\qed

\vspace{2mm}
(Proof of Theorem \ref{main1})

The condition $v_p(L^{\prime}(E,1)_{\alg} \cdot \Sub_{\al,\be}) \geq 1$ in Theorem \ref{main1} is equivalent to
\[
v_p(L^{\prime}(E,1)_{\alg}\cdot\Sub_p) \geq 1-v_p(\al-\be)
\] 
by Lemma \ref{subst}. We note that $v_p(\al-\be)$ is $0$ when $p\nmid a_p(E)$ and $1/2$  when $p\mid a_p(E)$. Since $v_p(L^{\prime}(E,1)_{\alg} \cdot\Sub_p) \in \Z$, the above inequality is equivalent to 
\[
v_p(L^{\prime}(E,1)_{\alg}\cdot\Sub_p) \geq 1. 
\]
By Lemma \ref{KEY}, this is also equivalent to 
\begin{align}\label{key ineq}
v_p(\#\Sha(E/\Q)^2\log_{\omega}(P)^2)\geq 4 \iff v_p(\#\Sha(E/\Q)\log_{\omega}(P))\geq 2 
\end{align}
 since $v_p(\#\Sha(E/\Q)^2\log_{\omega}(P)^2)$ is even. In the following, we divide the proof into two cases.

\vspace{2mm}
\noindent(Case 1 : $v_p(\#\Sha(E/\Q))\geq 1$)\ \ In this case, we have $\dim_{\F_p}(\Sha(E/\Q)[p])\geq 2$ by the existence of the Cassels-Tate pairing on $\Sha(E/\Q)$, and therefore $\dim_{\F_p}(\Sel(\Q, E[p]))\geq 2$. Then Proposition \ref{PS result} for $F=\Q$ implies that we have $r(E[p]) \neq 0$. Here the assumption (2) in Proposition \ref{PS result} holds by the assumption (4) in Theorem \ref{main1}.

\vspace{2mm}
\noindent(Case 2 : $v_p(\#\Sha(E/\Q))=0$)\ \ In this case, (\ref{key ineq}) implies that $v_p(\log_{\omega}(P)) \geq 2$. Since we assume $p\nmid f:=\#\widetilde{E}(\F_p)$, 
\[
v_p(\log_{\omega}(P)) = v_p(f\log_{\omega}(P)) = v_p(\log_{\omega}(fP)).
\]
By the exact sequence
\[
0 \to E_1(\Q_p) \to E(\Q_p) \to \widetilde{E}(\F_p) \to 0,
\]
the rational point $fP$ is in $(E(\Q)\cap E_1(\Q_p)) \setminus pE(\Q)$. Applying Proposition \ref{my result} for $F=\Q$ and $Q=fP$, we obtain $r(E[p]) \neq 0$. \hfill\qed

\begin{rem}\label{imp of PS}
By combining Lemma \ref{subst} and the equality (\ref{summary}), we obtain 
\[
v_p(\Sub_{\al, \be}) = v_p(\#\Sha(E/\Q)\log_{\omega}(P)^2)+v_p(\al-\be)-2.
\]
Here $v_p(\log_{\omega}(P)) \geq  1$ by the assumption $(4)$ in Theorem \ref{main1}. Therefore, we have $v_p(\Sub_{\al, \be})\geq 0$ and Theorem \ref{main1} refines the consequence (\ref{intromot}) of the result by Prasad and Shekhar \cite{PS21} under the assumptions in the theorem, as we noted in Remark \ref{introrefine}. 
\end{rem}

\subsection{Result on division fields over imaginary quadratic fields}
Let $E$ be an elliptic curve over $\Q$ which has good reduction at an odd prime number $p$. Let $K$ be an imaginary quadratic field which satisfies the three conditions (a), (b), and (c) at the beginning of Subsection \ref{antipL}.
\begin{thm}[Theorem \ref{intromain2} in Section \ref{introresults}]\label{main2}
We assume the following conditions: 
\begin{itemize}
\item[(1)] The elliptic curve $E$ has good reduction at $p$ and $\ra(E/K)=1$.
\item[(2)] The anticyclotomic Iwasawa main conjecture \cite[Conjecture 2.4]{CHKLL23} for $E$ and $p$ holds. 
\item[(3)] The prime number $p$ does not divide $\# E(\F_p)\cdot \mathrm{Tam}(E/\Q)$.
\item[(4)] As a $\Gal(K_E/K)$-module, $E[p]$ is irreducible, where $K_E:=K(E[p])$.
\end{itemize}
Then we have 
\[
v_p\left(L^{\prime}(E/K, 1)_{\alg} \cdot L^{\mathrm{BDP}}_{\mathfrak{p}}(0)\right) \geq 1 \hspace{2mm}\Rightarrow \hspace{2mm} r_K(E[p]) \neq 0,
\]
where the algebraic part $L^{\prime}(E/K, 1)_{\alg}$ of the derivative $L^{\prime}(E/K, 1)$ is defined in Section \ref{introimag}.
\end{thm}
\begin{rem}\label{an alg K}
By the work of Gross--Zagier and Kolyvagin, we have $\ral(E/K):=\mathrm{rank}_{\Z}(E(K))=1$ in Theorem \ref{main2} since $\ra(E/K)=1$.
\end{rem}

\proofname\ \ By the assumptions $(1), (2), (4)$ and Proposition \ref{antiIMCtopaBSD}, the $p$-part of Conjecture \ref{antipaBSD} now holds, which asserts that 
\begin{align*}
v_p(L_{\mathfrak{p}}^{\mathrm{BDP}}(0)) = v_p\left(\left(1-a_p(E)p^{-1}+p^{-1}\right)^2 \cdot \#\Sha(E/K)\cdot \log_{\om}(P)^2 \cdot \mathrm{Tam}(E/\Q)^2\right).
\end{align*}
Here $P \in E(K)$ denotes a generator of $E(K)_{/\mathrm{tors}}$. As in the proof of Theorem \ref{main1}, the $p$-adic valuations of the factors $\mathrm{Tam}(E/\Q)^2$ and $1-a_p(E)p^{-1}+p^{-1}$ are $0$ and $-1$, respectively by the assumption $(3)$. Thus we have 
\begin{align}\label{c of antipaBSD}
v_p(L_{\mathfrak{p}}^{\mathrm{BDP}}(0))=v_p(\#\Sha(E/K)\cdot \log_{\om}(P)^2)-2
\end{align}
On the other hand, by \cite[Proposition 4.3 (i)]{CHKLL23}, the $p$-part of Conjecture \ref{antipaBSD} is equivalent to the $p$-part of the BSD conjecture for $L(E/K, s)$, which asserts
\[
v_p(L^{\prime}(E/K, 1)_{\alg}) = v_p\left(\frac{\#\Sha(E/K)\cdot\mathrm{Tam}(E/\Q)^2}{\#E(K)^2_{\mathrm{tors}}}\right)=v_p(\#\Sha(E/K)).
\]  
Combining this equality with (\ref{c of antipaBSD}), we obtain
\[
v_p(L^{\prime}(E/K, 1)_{\alg}\cdot L_{\mathfrak{p}}^{\mathrm{BDP}}(0)) = v_p(\#\Sha(E/K)^2\cdot \log_{\om}(P)^2)-2.
\]
By this equality, if the condition $v_p(L^{\prime}(E/K, 1)_{\alg}\cdot L_{\mathfrak{p}}^{\mathrm{BDP}}(0))\geq 1$ in Theorem \ref{main2} holds, we have 
\[
v_p(\#\Sha(E/K)^2\cdot \log_{\om}(P)^2)\geq 4 \iff v_p(\#\Sha(E/K)\cdot \log_{\om}(P))\geq 2,
\]
since $v_p(\#\Sha(E/K)^2\cdot \log_{\om}(P)^2)$ is an even number. In the following, we divide the proof into two cases.
 
\vspace{2mm}
\noindent(Case 1 : $v_p(\#\Sha(E/K))\geq 1$)\ \ In this case, we have $\dim_{\F_p}(\Sha(E/K)[p])\geq 2$ by the existence of the Cassels-Tate pairing on $\Sha(E/K)$. We consider the following exact sequence
\[
0 \to E(K)/pE(K) \to \Sel(K, E[p]) \to \Sha(E/K)[p] \to 0.
\]
Since $\ral(E/K) =1$ by Remark \ref{an alg K}, we have $\dim_{\F_p}(E(K)/pE(K)) = 1$, which implies $\dim_{\F_p}(\Sel(K, E[p])) \geq 3$ by the above exact sequence. Using Proposition \ref{PS result} for $F=K$, we have $r_K(E[p])\neq 0$. Here the assumption (2) in Proposition \ref{PS result} holds by the assumption (3) in Theorem \ref{main2}.

\vspace{2mm}
\noindent(Case 2 : $v_p(\#\Sha(E/K))=0$)\ \ In this case, we have $v_p(\log_{\omega}(P)) \geq 2$. Since we assume $p\nmid f:=\#E(\F_p)$ in the condition $(2)$, 
\[
v_p(\log_{\omega}(P)) = v_p(f\log_{\omega}(P)) = v_p(\log_{\omega}(fP)).
\]
Then the rational point $fP$ is in $(E(K)\cap E_1(K_{\mathfrak{p}})) \setminus pE(K)$ by the same argument as in the proof of Theorem \ref{main1}. Since $K$ is an imaginary quadratic field in which $p$ splits, all the assumptions in Proposition \ref{my result} hold now. Thus applying Proposition \ref{my result} for $F=K$ and $Q=fP$, we obtain $r_K(E[p]) \neq 0$. \hfill\qed

\begin{example}\label{anti example}
We note that the value $L^{\mathrm{BDP}}_{\mathfrak{p}}(0)$ lies outside of the range of interpolation, and we cannot compute the value from the complex $L$-function of $E/K$. 

Let $p=31$ and $E$ be the elliptic curve defined by the equation $y^2+xy+y=x^3-x^2$ whose LMFDB \cite{lmfdb} label is 53.a1. The conductor of $E$ is $53$. The elliptic curve $E$ has good reduction at $31$ with $31\nmid a_{31}(E)$. 

Put $K:=\Q(\sqrt{-11})$. This imaginary quadratic field $K$ satisfies all the conditions (a), (b), and (c) at the beginning of Subsection \ref{antipL}.  
The conditions (1), (3)  and (4) in Theorem \ref{main2} can be checked numerically. Further, we can see that the mod $31$ Galois representation $\rho_{E, 31} : \Gal(K_E/K) \to \mathrm{Aut}_{\F_{31}}(E[31])$ is surjective. Then the condition (2) is confirmed by \cite[Theorem B.]{BCK21} under the condition (3) and the surjectivity of $\rho_{E, 31}$.

In this case, we can check that $v_{31}(L^{\prime}(E/K, 1)_{\alg}) =0$. As we described in the proof of Theorem \ref{main2}, we have
\[
v_{31}(L^{\prime}(E/K, 1)_{\alg}\cdot L_{\mathfrak{p}}^{\mathrm{BDP}}(0))= v_{31}(L_{\mathfrak{p}}^{\mathrm{BDP}}(0)) = v_{31}(\log_{\om}(P)^2)- 2,
\]
where $P$ is a generator of the torsion-free quotient $E(K)_{/\mathrm{tors}}$ of $E(K)$. We take $(0,0) \in E(K)$ as $P$. Then putting $\#\tilde{E}(\F_{31})\cdot (0,0) = 28 \cdot (0,0)$ as $(X, Y)$, we have 
\[
X =\frac{41\cdot 1811 \cdot 4019 \cdot 20047\cdot 511337\cdot 12164057487289}{2^2\cdot 3^4 \cdot \bm{31}^4\cdot 607^2\cdot 467258663^2},
\]
\[
Y=-\frac{5^2\cdot 41^2\cdot 61\cdot 239\cdot 4409\cdot 27329\cdot 32251\cdot 12164057487289^2}{2^3\cdot 3^6\cdot\bm{31^6}\cdot 607^3\cdot 467258663^3}.
\]
The above computation shows $v_{31}(\log_{\omega}(28 \cdot (0,0))) =2$, and we obtain $v_{31}(L^{\prime}(E/K, 1)_{\alg}\cdot L_{\mathfrak{p}}^{\mathrm{BDP}}(0)) = 2$. Then the non-vanishing of $r_K(E[31])$ follows from Theorem \ref{main2}.
\end{example}

\subsection*{Acknowledgement}
The author would like to thank Masato Kurihara for his continued support and valuable comments. The author would like to thank Takamichi Sano for his helpful comments on Perrin-Riou's work. The author also thanks Shinichi Kobayashi heartily for giving him some valuable comments on this work.

\bibliographystyle{habbrv.bst} \bibliography{padicLidealclassgroups}
\end{document}